\documentclass[MSNbibl,number,citesort,seceqn,dvips]{arxbj}
\usepackage{mathrsfs}
\aid{0}
\volume{19}
\issue{1}
\pubyear{2013}
\firstpage{44}
\lastpage{63}
\doi{10.3150/11-BEJ385}

\makeatletter

\def\bsuffix #1{#1}

\newtheorem{theorem}{Theorem}[section]
\newtheorem{lem}[theorem]{Lemma}

\newremark{rem}[theorem]{Remark}
\newremark{ex}[theorem]{Example}
\newproclaim{df}[theorem]{Definition}

\makeatother

\begin{document}
\begin{frontmatter}

\title{Consistent nonparametric Bayesian inference for
discretely observed scalar diffusions}
\runtitle{Consistent Bayesian inference for diffusions}

\begin{aug}
%%%% inicialai - be tarpu
\author[1]{\fnms{Frank} \snm{van der Meulen}\corref{}\thanksref{1}\ead[label=e1]{f.h.vandermeulen@tudelft.nl}}
\and
\author[2]{\fnms{Harry} \snm{van Zanten}\thanksref{2}\ead[label=e2]{j.h.v.zanten@tue.nl}}
\runauthor{F. van der Meulen and H. van Zanten} %% auto
\address[1]{Delft Institute of Applied Mathematics (DIAM),
Faculty of Electrical Engineering, Mathematics and Computer Science,
Delft University of Technology,
Mekelweg 4, 2628 CD Delft,
The Netherlands. \mbox{\printead{e1}}}
\address[2]{Department of Mathematics,
Eindhoven University of Technology,
P.O. Box 513,
5600 MB Eindhoven,
The~Netherlands. \printead{e2}}
\end{aug}

% HISTORY:
\received{\smonth{3} \syear{2010}}
\revised{\smonth{2} \syear{2011}}

% ABSTRACT
%
\begin{abstract}
We study Bayes procedures
for the problem of nonparametric drift estimation
for one-dimensional, ergodic diffusion models
from discrete-time, low-frequency data.
We give conditions for posterior consistency
and verify these conditions
for concrete priors,
including priors based on % integrated Brownian motion or
wavelet expansions.
\end{abstract}

% KEYWORDS
%

\end{frontmatter}

\section{Introduction}\label{sec1}

Consider the one-dimensional diffusion model
\begin{equation}\label{eqsde}
\mathrm{d}X_t = b(X_t) \,\mathrm{d}t + \mathrm{d}W_t, \qquad t \ge0,
\end{equation}
where $W$ is a standard Brownian motion and $b$ is an unknown drift
function belonging to a class
of functions $\mathscr{B}$.
We will make assumptions on $b$, stated precisely in the next section,
ensuring that (\ref{eqsde}) has a unique stationary solution $X$.
The aim is to make inference about $b$ on the
basis of discrete-time observations $X_0, X_\Delta, \ldots,
X_{n\Delta}$,
for some fixed sampling frequency $1/\Delta$.

Under appropriate conditions the solution to (\ref{eqsde}) is a
positively recurrent,
ergodic Markov process with a unique invariant probability distribution.
Moreover, under mild regularity conditions the process has transition
densities $p_b(t,x,y)$
relative to Lebesgue measure. In this case, we can employ a Bayes procedure
for making inference about the drift
function $b$. This involves putting a prior distribution $\Pi$ on the
set of drift functions $\mathscr{B}$ and
computing the posterior $\Pi(\cdot|X_0, X_\Delta, \ldots,
X_{n\Delta})$.
If the initial distribution is the invariant probability measure with
density $\pi_b$, the posterior
measure of a measurable set $B \subset\mathscr{B}$ is given by
\begin{eqnarray*}
\Pi(B|X_0, \ldots, X_{n\Delta}) =
\frac{\int_B \pi_b(X_0)\prod_{i=1}^n p_b(\Delta, X_{(i-1)\Delta},
X_{i\Delta}) \Pi(\mathrm{d}b)}
{\int_\mathscr{B}\pi_b(X_0)\prod_{i=1}^n p_b(\Delta,
X_{(i-1)\Delta},
X_{i\Delta}) \Pi(\mathrm{d}b)}.
\end{eqnarray*}
(We assume of course the necessary measurability to ensure that this is
well defined.)

This immediately reveals a practical complication, since the transition
densities of a diffusion
process can typically not be computed explicitly.
Several approaches have been
proposed in the literature to circumvent this problem.
These include for instance simulation-based methods for approximating
the transition densities,
or Y. A\"it-Sahalia's closed-form expansions, cf. for example, Jensen and Poulsen \cite{Jensen} for an overview.
A method that has been
proven to be particularly useful for dealing with Bayes procedures is
to view the continuous segments of the diffusion process
between the observations as missing data and to employ a Gibbs sampling scheme.
Practically this involves repeatedly simulating diffusion bridges to
generate the missing
data and drawing from the posterior distribution of $b$ given the augmented,
continuous data $(X_t\dvt t\in[0,n\Delta])$.
Several schemes have been devised to simulate the diffusion bridges,
see, for example, Elerian \textit{et~al.} \cite{Elerian01}, Eraker \cite{Eraker01},
Roberts and Stramer \cite{Roberts01}, Beskos \textit{et~al.} \cite{Bes06},
Golightly and Wilkinson \cite{Golightly} and Chib \textit{et al.} \cite{Chib}. Drawing
from the continuous-data posterior can be
done by more conventional methods, because contrary to the
discrete-observations likelihood,
the continuous-data likelihood has a known closed form expression given
by Girsanov's theorem.

For parametric models, where the drift function is known up to a
Euclidean parameter $\theta$
that has to be estimated, the outlined approach has been shown to
provide an effective method
for dealing with discretely observed diffusions.
The approach is however not essentially limited to a parametric setup.
The methodology has great potential to be developed into a practically
feasible methodology
in nonparametric settings as well.
It is however very well known that in Bayesian nonparametrics the
choice of
the prior distribution is crucial and posterior consistency is not
automatically guaranteed
(e.g., Diaconis and Freedman \cite{Diac86}). This motivates the study of posterior
consistency for discretely observed
diffusions carried out in this paper.

In the i.i.d.-setting, sufficient conditions for posterior consistency
were first obtained by Schwartz \cite{Schwartz}. See also Barron \textit{et al.} \cite{barron}, Ghosal and van~der Vaart \cite{vdv}
and Shen and Wasserman \cite{shenwas}. Here we consider discrete observations from
a diffusion model (\ref{eqsde}), which constitute a
Markov chain. A number of recent papers have investigated the
problem of posterior consistency or convergence rates for Markov data,
cf.
for example, Ghosal and van~der Vaart \cite{Gho07}, Ghosal and Tang \cite{TanGhos}, Tang and Ghosal \cite{Tan07}.
The results in these papers do however
not immediately lead to practically useful results for our setting. The
problem lies again in
the fact that in our case, the transition densities of the model are typically
not analytically tractable.
Since the conditions for consistency given for instance by  Tang and Ghosal \cite{Tan07}
involve the transition densities, they cannot be readily used
to verify consistency for a given prior in our discretely observed
diffusion model.
The aim in this paper is to formulate conditions involving only the
coefficients appearing in the stochastic differential equation (\ref{eqsde}).
We achieve this by adapting the results of Tang and Ghosal \cite{Tan07} to the present
setting. Basically, we need two assumptions. Firstly, if $\mu_0$
denotes the true invariant probability measure of the process $X$,
we require that the prior puts positive mass on balls $\{b\in{\mathscr
B} \dvt  \|b-b_0\|_{2,\mu_0} < \varepsilon\}$ for each $\varepsilon>0$
($\|\cdot\|
_{2,\mu_0}$ denotes a weighted $L^2$-norm and $b_0$ denotes the true
drift). This is a natural condition, since if the prior excludes the
true drift, consistency can never be obtained.
Secondly, we need an equicontinuity assumption (Definition~\ref{dflocunifequicont}), which limits the size, or rather the
complexity, of the set of drift functions.
Under these assumptions, we obtain posterior consistency (Theorem \ref{theoremcon}):
the posterior measure of appropriately defined weak neighborhoods of
the true drift function $b_0$ converges to $1$
almost surely, as the number of observations $n$ tends to infinity.
This is the main result of the paper.

Ghosal and van~der Vaart \cite{Gho07} give conditions from which the posterior rate of
convergence for Markov chain data can be calculated. These conditions
are a combination of a prior mass condition and a testing condition.
This testing condition requires that one can test the true drift
function against balls of alternatives with exponentially decaying
error probabilities. Such tests are not easily constructed in the
present setup. Appropriate tests for Markov chains have been shown to
exist under certain (lower) bounds on the transition probabilities
(e.g., \cite{Birge}).
In our setup such bounds do however not seem to be valid in general. An
interesting line of future research would be to extend or adapt the
available testing results for Markov chains to the setting of
discretely observed diffusions.
This may not only give posterior consistency results in a stronger
topology, but may pave the way for obtaining posterior rates of
convergence as well.
In the present paper, we completely avoid the construction of tests.
Instead we employ martingale arguments
in a similar fashion as Tang and Ghosal \cite{Tan07}, who adapted the approach of
Walker~\cite{Walker} to the Markov chain setting.

%
%The problem of nonparametric drift estimation for diffusions based on
%discrete time, low-frequency observations has also been addressed in a
%non-Bayesian setup. \cite{Gobet} considered diffusions having a
%generator
%with countable spectrum, such as diffusions on a bounded interval with
%reflecting boundary conditions.
%Under these assumptions, they constructed estimators for the drift and
%diffusion functions that attain optimal
%convergence rates in the low-frequency observations setting.
%Another basic idea to make inference about the drift function is to
%consider
%the pseudo-observations $Y_1,\ldots,Y_n$, where $Y_i=(X_{i
%drift function $b$ can be obtained by regressing $Y_i$ on $X_{i
%kernel estimators, smoothing splines or penalized splines. Of course,
%this approach crucially depends on a sufficiently high sampling
%frequency. In any case, already in the parametric case, the resulting
%methods are inconsistent
%if the sampling frequency $1/\Delta$ is fixed.
%In the present setup, kernel estimators (Nadaraya-Watson estimator)
%have been considered by e.g.\ \cite{Stanton97} and \cite{Jiang97}. To
%remove boundary effects, \cite{Fan03} proposed local linear
%estimation. \cite{Comte07} consider penalized least squares estimation.
%For an overview of various nonparametric estimation methods for
%diffusions we refer to \cite{Fan05} and \cite{Cai03}.

The remainder of the paper is organized as follows.
In Section \ref{secsetup}, preliminaries on the statistical model
and Bayes procedure are outlined. The main consistency result of this
paper is formulated in Section \ref{seccon}. Examples of priors that
satisfy the requirements for consistency are given in Section~\ref{secexamples}. The paper ends with a proof of the main result and some
concluding remarks. The \hyperref[app]{Appendix} contains a technical lemma.

\subsection{Notation}

$\|g\|_{p,\nu} = (\int|g|^p \,\mathrm{d}\nu)^{1/p}$: $L_p$-norm relative to
the measure $\nu$.

\noindent$L^2(\mu)$: space of square integrable functions with
respect to measure $\mu$.

\noindent$C(A)$, $BC(A)$: space of continuous functions, space of
bounded continuous functions defined on $A \subseteq\mathbb{R}$.

\noindent
$C^s(\mathbb{R})$, for $s \in(0,1)$: space of $s$-H\"older functions, that
is,
\begin{eqnarray}\label{eqholder}
C^s(\mathbb{R}) =\biggl\{ f \in C(\mathbb{R}) \dvt \|f\|_s = \sup_{x, h}
\frac
{|f(x+h)-f(x)|}{|h|^s} < \infty\biggr\}.
\end{eqnarray}

%and for $s=n+s'$ ($n \in\{1,2\ldots\}$, $s' \in(0,1)$) by
%$$ C^s(\RR)= \{ f \in C^n(\RR) ~ : ~ \frac{d^n}{d x^n} f \in
%C^{s'}\}. $$
%Note that these definitions exclude $C^s(\RR)$ for $s$ equal to a
%positive integer. For reasons %explained in remark

%(0,1)$ and radius $L > 0$, i.e.\
%$C^\alpha_L[-m,m] = \{f: [-m,m] \to\RR: \sup_{x,y \in[-m,m]}
%|f(x)-f(y)|/|x-y|^\alpha\le L\}$.

\noindent${\mathscr L}(X)$: law of a random variable $X$.

\noindent$\mathbb{P}^b_\mu$: law that the solution of the SDE (\ref{eqsde}), with ${\mathscr L}(X_0)=\mu$, generated
on the canonical path space $C(\mathbb{R}_+)$.

\noindent$\mu_ b$: invariant measure.

\noindent$\mathbb{P}_b$: short-hand notation for $\mathbb{P}^b_{\mu_b}$.

\noindent$\mathbb{P}_x^b$: short-hand notation for $\mathbb
{P}^b_{\delta_x}$,
where $\delta_x$ denotes Dirac measure at $x\in\mathbb{R}$.

\noindent$\mu_0$: short-hand notation for $\mu_{b_0}$.

\noindent$\pi_b$: density of invariant measure.

\noindent$(P^b_t)_{t \ge0}$: transition semigroup associated with the
diffusion.

\noindent$p_b(t,x,y)$: transition density.

\section{Setup}
\label{secsetup}

\subsection{Description of the diffusion model}

In this section, we give a precise description of the diffusion model
that we consider.
Let $\mathscr{B}\subset C(\mathbb{R})$ be a collection of continuous
functions on
$\mathbb{R}$.
For $b \in\mathscr{B}$ and a fixed number $c\in\mathbb{R}$, let the
function $s_b\dvtx
\mathbb{R}\to\mathbb{R}$ be defined by
\[
s_b(x) = \int_{c}^x\exp\biggl(-2\int_{c}^y {b(z)} \,\mathrm{d}z\biggr) \,\mathrm{d}y.
\]
We assume that
\[
\lim_{x\downarrow-\infty} s_b(x) = -\infty, \qquad\lim_{x\uparrow
\infty} s_b(x) = \infty
\]
for all $b \in\mathscr{B}$. The finiteness (or nonfiniteness) of
these limits
does not depend on the choice of $c$ (see page 339 in Karatzas and Shreve \cite{karatzas}).
It is classical that under these assumptions, we have that for every $x
\in\mathbb{R}$ and $b \in\mathscr{B}$,
the SDE
\[
\mathrm{d}X_t = b(X_t) \,\mathrm{d}t + \mathrm{d}W_t, \qquad X_0 = x,
\]
has a unique weak solution. Let $\mathbb{P}^b_x$ denote the law that this
solution generates
on the canonical path space $C(\mathbb{R}_+)$. Then in the commonly used
terminology of
It{\^o} and McKean \cite{Ito1} or Kallenberg \cite{Kal1}, Chapter 23, the collection of laws
$(\mathbb{P}
_x^b\dvt x \in\mathbb{R})$
constitutes a canonical, recurrent diffusion on the real line. In other
words, for $X$
the canonical process on $\Omega= C(\mathbb{R}_+)$ defined by
$X_t(\omega)
= \omega(t)$
we have the following:
\begin{longlist}[(iii)]
\item[(i)] Under $\mathbb{P}^b_x$ the process $X$ starts in $x$, that is,\
$\mathbb{P}
^b_x(X_0 = x) = 1$
for all $x \in\mathbb{R}$.

\item[(ii)] The process $X$ is strong Markov.

\item[(iii)] For all $x \in\mathbb{R}$, the process $X$ is recurrent under
$\mathbb{P}
^b_x$.
\end{longlist}
For a probability measure $\mu$ on $\mathbb{R}$ we define, as usual,
$\mathbb{P}
^b_\mu(B) = \int\mathbb{P}^b_x(B) \mu(\mathrm{d}x)$ for a measurable set $B$.
Then under $\mathbb{P}^b_\mu$ the law of $X_0$ equals $\mu$ and $X$
is the
weak solution of
\[
\mathrm{d}X_t = b(X_t) \,\mathrm{d}t + \mathrm{d}W_t, \qquad\mathscr{L}(X_0) = \mu.
\]

As the notation suggests, $s_b$ is the scale function of the diffusion.
The speed measure is denoted by $m_b$.
In the present setting, it is the Borel measure on $\mathbb{R}$ given by
\begin{equation}\label{eqspeed}
m_b(\mathrm{d}x) = \exp \biggl(2\int_{0}^x{b(z)} \,\mathrm{d}z\biggr) \,\mathrm{d}x.
\end{equation}
We assume that the speed measure is finite, that is, $m_b(\mathbb{R}) <
\infty
$. This ensures that the diffusion
is positively recurrent and ergodic in the sense that for all $x \in
\mathbb{R}$,
\begin{equation}\label{eqerg}
X_t \stackrel{\mathbb{P}^b_x}{\Longrightarrow} \mu_b
\end{equation}
as $t \to\infty$, where $\mu_b = m_b/m_b(\mathbb{R})$ is the normalized
speed measure
(cf., e.g., Kallenberg \cite{Kal1}, Theorem 23.15). We will write $\mu_{0}=\mu_{b_0}$.
The measure $\mu_b$ is the unique invariant probability measure of the
diffusion.
In particular, the process $X$ is stationary under $\mathbb{P}^b_{\mu_{b}}$.
It is easily
verified that under our conditions,
$\mu_b$ has a continuously differentiable Lebesgue density $\pi_b$.
Moreover, it follows
from (\ref{eqspeed}) that we have the
relation
\begin{equation}\label{eqdriftrelation}
b = \frac{\pi'_b}{2\pi_b}.
\end{equation}

We denote the transition semigroup associated to the diffusion by
$(P^b_t)_{t \ge0}$.
In other words, for a bounded measurable function $f$ on $\mathbb{R}$
and $x
\in\mathbb{R}$ we have
$P_t^b f (x) = {\mathbb{E}}^b_xf(X_t)$, where ${\mathbb{E}}^b_x$ is
the expectation associated
to $\mathbb{P}^b_x$.
The operator $P^b_t$ maps the space $BC(\mathbb{R})$ of bounded, continuous
functions on $\mathbb{R}$ into itself
(see, e.g., (the proof of) Theorem 23.13 of Kallenberg \cite{Kal1}, or Rogers and Williams \cite{Rog1}, Proposition V.50.1).
A regular diffusion as we are considering is known to have
positive transition densities with respect to its speed measure, cf.,
for example, It{\^o} and McKean \cite{Ito1}, Section 4.11.
Since the speed measure has a positive Lebesgue density under our
assumptions, we have in fact
the existence of transition densities $p_b\dvt (0,\infty) \times\mathbb{R}
\times\mathbb{R}\to(0,\infty)$ such that
for all bounded, measurable functions $f$, $x \in\mathbb{R}$ and $t > 0$,
\[
P^b_t f(x) = \int_\mathbb{R}p_b(t,x,y)f(y) \,\mathrm{d}y.
\]

For more background on the theory of one-dimensional diffusions and relevant
references to the literature, see, for instance, Borodin and Salminen \cite{Bor02}.

\subsection{Statistical model and Bayes procedure}

Consider the setting described in the preceding section,
that is, we have a collection $\mathscr{B}\subset C(\mathbb{R})$
such that
every $b \in\mathscr{B}$ determines an SDE that generates an ergodic
diffusion on $\mathbb{R}$.
For $b \in\mathscr{B}$, let $\mathbb{P}_b$ be defined by
$\mathbb{P}_b = \mathbb{P}_{\mu_b}^b$. In other words, under
$\mathbb{P}_b$ the
canonical process $X$ on $C(\mathbb{R}_+)$ is the unique stationary
solution of
the SDE
\[
\mathrm{d}X_t = b(X_t) \,\mathrm{d}t + \mathrm{d}W_t.
\]
We assume that for some fixed $\Delta> 0$ and a natural number $n$, we
have $n+1$
observations $X_0, X_\Delta, \ldots, X_{n\Delta}$ from $X$ under
$\mathbb{P}_{b_0}$,
for some
``true'' drift function $b_0 \in\mathscr{B}$.
The aim is to infer the drift function $b_0$ from these data.

In our Bayesian approach, we assume that the model $\mathscr{B}$ is a
measurable subset of $C(\mathbb{R})$ and
we put a prior distribution $\Pi$ on it.
Next, we consider the posterior distribution
$\Pi(\cdot|X_0, \ldots, X_{n\Delta})$ on $\mathscr{B}$, which is
given by
\[
\Pi(B|X_0, \ldots, X_{n\Delta}) =
\frac{\int_B \pi_b(X_0)\prod_{i=1}^n p_b(\Delta, X_{(i-1)\Delta},
X_{i\Delta}) \Pi(\mathrm{d}b)}
{\int_\mathscr{B}\pi_b(X_0)\prod_{i=1}^n p_b(\Delta,
X_{(i-1)\Delta},
X_{i\Delta}) \Pi(\mathrm{d}b)}.
\]
In the next section, we provide sufficient conditions under which the posterior
asymptotically concentrates its mass around the true drift function
$b_0$ as $n \to\infty$.

\section{Consistency}
\label{seccon}

We are interested in conditions under which the posterior
asymptotically concentrates its mass around the true drift function $b_0$.
More precisely, we want that under $\mathbb{P}_{b_0}$ the posterior
mass concentrates on arbitrarily small neighborhoods of $b_0$.
To ensure that neighborhoods of points $b \not= b_0$ do not receive
posterior mass in the limit, the topology we use to define the neighborhoods
should have some separation properties, it should for instance be Hausdorff.

We define a weak topology on $\mathscr{B}$ through the transition operators
$P^b_\Delta$
(see Section \ref{secsetup}).
This is justified by the following lemma, which states that identifying
the drift parameter $b$
is in our setting equivalent to identifying $P^b_\Delta$.

\begin{lem}\label{lemuniquedrift}
If $P^b_t = P^{b'}_t$ for some $t > 0$, then $b = b'$.
\end{lem}

\begin{pf}
Fix an $x \in\mathbb{R}$ and $b \in\mathscr{B}$.
By the semigroup property, the law of $X_{nt}$ under $\mathbb{P}^b_x$ is
determined by $P^b_t$.
Indeed, for $f$ a bounded measurable function and $n$ a natural number
we have
\[
{\mathbb{E}}^b_x f(X_{nt}) = (P^b_t)^nf(x).
\]
On the other hand, ergodicity implies that the law of $X_{nt}$ under
$\mathbb{P}^b_x$
converges weakly to the invariant distribution $\mu_b$, cf. (\ref{eqerg}).
It follows that $P^b_t$ completely determines $\mu_b$.
By (\ref{eqdriftrelation}), $\mu_b$ completely determines $b$ under
our assumptions.
\end{pf}

Now let $\nu$ be a finite Borel measure on the state space $\mathbb{R}$.
For $b \in\mathscr{B}$, $f \in BC(\mathbb{R})$ and $\varepsilon>
0$, let
\[
U^{b}_{f,\varepsilon} = \{b' \in\mathscr{B}\dvt \|P^{b'}_\Delta f -
P^{b}_\Delta f\|
_{1, \nu} < \varepsilon\}.
\]
Consider the topology on $\mathscr{B}$ that is determined by the requirement
that for $b \in\mathscr{B}$,
the collection of sets
\[
\{U^b_{f, \varepsilon}\dvt f \in BC(\mathbb{R}), \varepsilon> 0\}
\]
forms a sub-base for the neighborhood system at $b$. By definition,
this means that
any open neighborhood of $b \in\mathscr{B}$ is a union of finite
intersections of the form
$U^b_{f_1, \varepsilon_1} \cap\cdots\cap U^b_{f_m, \varepsilon_m}$.

Although the topology is defined in a rather indirect fashion, it
has the desired Hausdorff property, that is,\ different points in
$\mathscr{B}
$ can be separated
by disjoint open sets.

\begin{lem}
If $\nu$ assigns positive mass to all nonempty open intervals, then
the topology on $\mathscr{B}$ is Hausdorff.
\end{lem}

\begin{pf}
Consider two functions $b \not= b'$ in $\mathscr{B}$.
By Lemma \ref{lemuniquedrift}, we have $P^b_\Delta\not=
P^{b'}_\Delta$ and
hence there exists an $f \in BC(\mathbb{R})$ and an $x \in\mathbb
{R}$ such that
$P^b_\Delta f(x) \not= P^{b'}_\Delta f(x)$.
By continuity there exists in fact a nonempty open interval $J \subset
\mathbb{R}$ where the functions
$P^b_\Delta f$ and $P^{b'}_\Delta f$ are different. By the assumption
on $\nu$,
it follows that for some $\varepsilon> 0$,
\[
\|P^b_\Delta f - P^{b'}_\Delta f\|_{1, \nu} > \varepsilon.
\]
This implies that the neighborhoods $U^b_{f,\varepsilon/2}$ and
$U^{b'}_{f,\varepsilon/2}$
are disjoint.
\end{pf}

An alternative point of view on the topology that we use is obtained
by considering the high-frequency limit $\Delta\to0$.
Let $A_b$ be the generator of $X$ under $\mathbb{P}_b$, that is,\
$A_b f = bf' + f''/2$ for a $C^2$-function $f$. Then for small $\Delta$,
\[
P^{b_1}_\Delta f - P^{b_2}_\Delta f\approx\Delta(A_{b_1}f - A_{b_2} f)
= \Delta(b_1-b_2)f'.
\]
It follows that
for small $\Delta$, the constructed topology is close to the topology
induced by the $L^1(\nu)$-norm
on
the set of drift functions $\mathscr{B}$.

Having specified the topology, we can define {\em weak posterior
consistency}, or just {\em consistency}.
\begin{df}\label{defconsistency}
We have weak posterior consistency if for every open neighborhood
$U_{b_0}$ of $b_0$, it holds that
\[
\Pi(b \notin U_{b_0} |X_0, X_{\Delta}, \ldots, X_{n\Delta})
\to0
\qquad\mbox{$\mathbb{P}_{b_0}$-a.s.}
\]
as $n \to\infty$. Note that the word ``weak'' refers to the topology,
not to
the mode of stochastic convergence.
\end{df}

Theorem \ref{theoremcon} below is the main result of this section. It
needs the following definition.

\begin{df}\label{dflocunifequicont}
We call a collection $\mathscr{F}$ of real-valued functions on the
real line
\textup{locally uniformly equicontinuous} if for every $\varepsilon> 0$ and
every compact $K \subset\mathbb{R}$, there exists
a $\delta> 0$ such that
\[
\sup_{f \in\mathscr{F}}\mathop{\sup_{x,y \in K }}_{|x-y| < \delta}
|f(x)-f(y)| <
\varepsilon.
\]
\end{df}

%Note for instance that a class of functions $\FF$ is locally uniform
%equicontinuous if for every
%compact interval $K \subset\RR$, the class $\FF|_K$ consisting of the
%restrictions of the elements of
%$\FF$ to the interval $K$ is contained in a H\"older ball $C^
%(which may depend on $K$).
In Section \ref{secexamples} we give examples of locally uniformly
equicontinuous collections of functions.

\begin{theorem}\label{theoremcon}
Suppose we have discrete-time data from the stationary solution to the
stochastic differential equation
\begin{eqnarray*}
\mathrm{d}X_t = b(X_t) \,\mathrm{d}t + \mathrm{d}W_t, \quad t \ge0.
\end{eqnarray*}
Denote the invariant measure of the diffusion with drift $b_0$ by $\mu_0$.
%Consider the setup described in Section \ref{secsetup}.
Let $\Pi$ be a prior on the set of drift functions $\mathscr{B}$ and suppose
that $\mathscr{B}$ is locally uniformly equicontinuous and
$\sup_{b \in\mathscr{B}} \|b\|_\infty< \infty$.
Then if
\begin{equation}\label{eqkl}
\Pi(b \in\mathscr{B}\dvt \|{b- b_0}\|_{2,\mu_{_0}} < \varepsilon
) > 0
\qquad\mbox{for all $\varepsilon> 0$},
\end{equation}
we have weak consistency (as in Definition \ref{defconsistency}).
\end{theorem}

In Bayesian practice, a model set $\mathscr{B}$ is typically not
specified explicitly.
Usually some prior $\Pi$ is simply chosen and the procedure is carried
out. From this perspective,
the theorem states that if the chosen prior gives mass $1$ to a set of functions
that is uniformly bounded and locally uniformly equicontinuous, then we
have weak consistency
for every true $b_0$ in the $L^2(\mu_0)$-support of the prior.

Prior mass conditions like (\ref{eqkl}) are standard in results on
posterior consistency. Intuitively,
it is reasonable that if we want the posterior to concentrate around
$b_0$ asymptotically, the prior
should put sufficient mass near $b_0$ too. The uniform boundedness and
equicontinuity conditions
limit the size of the support of the prior, which is reasonable as
well. The conditions are
somewhat restrictive, but due to technical reasons cannot be avoided in
our approach.
In settings where consistency can be derived using testing arguments,
boundedness and equicontinuity
conditions can typically be relaxed, and only need to be valid on
certain subsets $\mathscr{B}_n$ of
the support $\mathscr{B}$ of the prior with increasing prior probability.
However, since we do not have the appropriate
tests available in this case, we cannot follow such an approach
unfortunately. On the other hand, computational
approaches like the one of Beskos \textit{et al.} \cite{Bes06} require in fact that both $b$
and its derivative $b'$ are uniformly bounded,
which is more restrictive than the conditions of our consistency theorem.

The proof of the theorem is deferred to Section \ref{secproof}.
In the next section, we first consider a number of concrete priors
for which the assumptions of the theorem are verified.

%-------------------------------------

\section{Examples of concrete priors}
\label{secexamples}

The following example is perhaps of little practical relevance, but it
shows already
that there is an abundance of priors available that yield posterior consistency.

\begin{ex}[(Discrete net priors)]
Let the collection of drift functions $\mathscr{B}$ satisfy the requirements
of Theorem \ref{theoremcon}. That is, $\mathscr{B}$ is locally uniformly
equicontinuous and $\sup_{b\in\mathscr{B}} \|b\|_\infty<\infty$.
To construct the prior choose
two probability distributions $(p_n)$ and $(q_n)$ on the positive
integers such that $p_n, q_n > 0$ for
$n$ large enough, and a decreasing sequence of positive numbers
$\varepsilon
_n \downarrow0$.
For $m \ge1$, let $\mathscr{B}_m = \{b|_{[-m, m]}\dvt b \in\mathscr
{B}\}$ be the set of
restrictions of functions in $\mathscr{B}$ to the interval $[-m,m]$.
The functions in $\mathscr{B}_m$ are uniformly equicontinuous and
hence, by
the Arzel\`{a}--Ascoli theorem,
$\mathscr{B}_m$ is totally bounded for the uniform norm.
For every $n$, we fix a finite $\varepsilon_n$-net $\mathscr{B}_{m,
\varepsilon_n}$ for
$\mathscr{B}_m$,
that is, $\mathscr{B}_{m, \varepsilon_n}$ is a finite set such that
every element
of $\mathscr{B}_m$ is within uniform
distance $\varepsilon_n$ of some element of $\mathscr{B}_{m,
\varepsilon_n}$.
We extend every function in the net to the whole real line by setting
it equal to $1$
on $(-\infty, -m-1]$ and to $-1$ on $[m+1, \infty)$, and
interpolating linearly in the intervals $[-m-1,-m]$
and $[m, m+1]$.
A draw $b$ from the prior $\Pi$ is now generated as follows:
\begin{longlist}[(iii)]
\item[(i)] draw $m$ from the probability distribution $(p_m)$,
\item[(ii)] draw $n$ from the probability distribution $(q_n)$,
\item[(iii)] draw $b$ uniformly from $\mathscr{B}_{m, \varepsilon_n}$.
\end{longlist}
In other words,
if $\mathscr{B}_{m, \varepsilon_n} = \{b^{m, n}_1, \ldots, b^{m,
n}_{k_{m,n}}\}$, then
\[
\Pi= \sum_{m=1}^{\infty} \sum_{n=1}^\infty\sum_{k=1}^{k_{m,n}}
\frac{ p_mq_n}{k_{m,n}} \delta_{b^{m,n}_k}.
\]

By construction, $\Pi$ assigns mass $1$ to a countable set of functions
that is uniformly bounded and locally uniformly continuous.
Now consider $b_0 \in\mathscr{B}$ and $\varepsilon>0$. We show that
condition (\ref{eqkl}) is satisfied.
For every $b \in\mathscr{B}$ and $m \in\mathbb{N}$, we have
\begin{eqnarray*}
\|b-b_0\|_{2,\mu_0}^2 & =& \int_{|x| \le m} \bigl(b(x)-b_0(x)\bigr)^2 \,\mathrm{d} \mu
_0(x) + \int_{|x|> m} \bigl(b(x)-b_0(x)\bigr)^2\, \mathrm{d} \mu_0(x)
\\
& \le&\|b-b_0\|^2_{m,\infty} + 2 \int_{|x|> m} \bigl(b^2(x)+ b_0^2(x)\bigr) \,\mathrm{d}
\mu_0(x)
\\
& \le&\|b-b_0\|^2_{m,\infty} + C \mu_0(|x| > m),
\end{eqnarray*}
where $\|\cdot\|_{m,\infty}$ denotes the uniform norm on $[-m,m]$ and
$C=2 (1+ \sup_{b\in\mathscr{B}} \|b\|_\infty^2) $.
Hence, for $m \in\mathbb{N}$ so large that $C \mu_0(|x| > m) <
\varepsilon^2$, it
holds that
\[
\Pi(b\dvt \|b-b_0\|_{2,\mu_0}^2 < 2\varepsilon^2 )
\ge\Pi(b\dvt \|b-b_0\|^2_{m,\infty} < \varepsilon^2).
\]
For $n \in\mathbb{N}$ such that $\varepsilon_n < \varepsilon$ and
$q_n > 0$, we have, by
construction,
\[
\Pi(b\dvt \|b-b_0\|^2_{m,\infty} < \varepsilon^2) \ge\Pi(b\dvt \|b-b_0\|
_{m,\infty} < \varepsilon_n) \ge\frac{p_mq_n}{k_{m,n}} > 0.
\]
This shows that condition (\ref{eqkl}) holds, and hence we have
posterior consistency for this class of priors.
\end{ex}

If $\mathscr{B}\subset C^s(\mathbb{R})$ for some $s \in(0,1)$ and
$\sup_{b \in\mathscr{B}} \|b\|_s < \infty$ (see (\ref{eqholder})), then
clearly $\mathscr{B}$ satisfies the
equicontinuity condition of Definition \ref{dflocunifequicont}. In
the following example,
we use wavelet expansions to construct a consistent prior on drift
functions which belong to such a class of H\"older functions.

\begin{ex}[(Wavelets)]\label{exwaveltes}
Suppose $\{ \varphi_{k}, \psi_{j,k}\}_{k\in\mathbb{Z}, j\ge0}$ is an
orthonormal wavelet basis, so that functions $f \in L^2(\mathbb{R})$
can be
represented as
\[
f = \sum_{k\in\mathbb{Z}} \langle f, \varphi_{k}\rangle\varphi
_{k} + \sum_{k
\in\mathbb{Z}} \sum_{j\ge0} \langle f, \psi_{j,k}\rangle\psi_{j,k}
\]
(the convergence being in $L^2(\mathbb{R})$). The functions $\psi
_{j,k}$ are
obtained from the mother wavelet function by translation and scaling:
$\psi_{j,k}(\cdot)= 2^{j/2} \psi(2^j \cdot-k)$. Similarly, the
$\varphi_{k}$ are obtained from the father wavelet $\varphi$ (also called
scaling function) by translation: $\varphi_{k}(\cdot)=\varphi(\cdot-k)$.

%It is well known that under appropriate smoothness conditions on $
%the smoothness of the function $f$. Assume $\psi$ is continuously
%differentiable %and that there exist constants $\eps>0$, $\gamma>0$,
%$C_0 <\infty$ and $C_1<\infty$ such that

It is well known that under appropriate smoothness conditions on $\psi
$, the rate of decay of the wavelet coefficients characterizes the
smoothness of the function $f$. Assume $\psi$ is continuously
differentiable and has compact support.
Then $f \in C^s \cap L^2(\mathbb{R})$ if and only if $\|f\|_\infty<
\infty$ and
\begin{itemize}
\item$|\langle f, \varphi_{k}\rangle| \le C_f$ for all $k\in\mathbb{Z}$,
\item$|\langle f, \psi_{j,k}\rangle| \le C_f 2^{-j(s+1/2)}$ for all
$j\ge0$ and $k\in\mathbb{Z}$.
\end{itemize}
Moreover, $C_f$ can be taken as the product of the H\"older norm of $f$
and a constant (that does not depend of $f$). For a proof, we refer to
Section 6.7 in Hern{\'a}ndez and Weiss \cite{hernandez}, see also Daubechies \cite{daub}, Section 9.2.
This characterization implies that for $s \in(0,1)$ and $L > 0$, the collection
\begin{eqnarray*}
 {\mathscr F}_{s,L} &:=& \biggl\{ f \in L^2(\mathbb{R}) \dvt f = \sum
_{k\in\mathbb{Z}}
a_{k} \varphi_{k} + \sum_{k \in\mathbb{Z}} \sum_{j\ge0} b_{j,k}
\psi
_{j,k},
\\
&&{}\hspace*{56pt}
{\sup_k} | a_{k} | + \sup_j \sup_k 2^{j(s+1/2)} |b_{j, k}| \le L\biggr\}
\end{eqnarray*}
consists of $s$-H\"older continuous functions with uniformly bounded H\"
older norms.

In addition to the smoothness condition on $\psi$, we assume that the
scaling function $\varphi$ is bounded and compactly supported. This
implies that the function $\theta_\varphi(x)=\sum_k | \varphi
(x-k)|$ is such
that $\mbox{ess sup}_x \theta_\varphi(x) < \infty$. This is a
localization condition that is referred to as Condition $(\theta)$ in
H{\"a}rdle \textit{et al}. \cite{hardle} (page 77). By inequalities (9.34) and (9.35) on page 114
in H{\"a}rdle \textit{et al.} \cite{hardle}, Condition $(\theta)$ implies that the supremum norm of
$\sum_k a_{k} \varphi_{0,k}$ is equivalent to the $\sup$-norm on the
sequence $\{a_{k}\}_{k}$. In addition, the supremum norm of $\sum
_{j\ge0} \sum_k b_{j,k} \psi_{j,k}$ is equivalent to the $\|\cdot\|
_w$-norm of the doubly indexed sequence $b=\{b_{j,k}\}_{j\ge0, k}$, where
\[
\|b\|_w = \sum_{j\ge0} 2^{j/2} \sup_k |b_{j,k}|.
\]
It follows in particular that the uniform norm of the functions in
$\mathscr{F}
_{s,L}$ is uniformly bounded.

To construct a prior on drift functions that is consistent for all true
drift functions $b_0 \in\mathscr{F}_{s, L}$ we first construct an
auxiliary prior $\Pi'$ on the whole class $\mathscr{F}_{s, L}$ (which does
not only charge drift functions of ergodic diffusions).
Let $J$ be a discrete random variable, supported on $\mathbb{N}_0= \{0,
1,\ldots\}$ and let $U_{j,k}, V_k$, for $j \in\mathbb{N}_0,
k \in\mathbb{Z}$, be independent random variables, independent of
$J$, from
a distribution with
a strictly positive, continuous density on its
support $[-L, L]$.
Define the prior $\Pi'$ as the law of the random function
\[
x \mapsto\sum_{k\in\mathbb{Z}} V_k \varphi_{k} + \sum_{j=0}^J
\sum_{k\in
\mathbb{Z}} \eta_j U_{j,k} \psi_{j,k}
\]
on $\mathbb{R}$, where $\eta_{j} = 2^{-j(s+1/2)}$. To arrive at a
prior on
drift functions of ergodic diffusions we
proceed as in the preceding example. We choose a probability
distribution $(p_m)$ on $\mathbb{N}$, with $p_m > 0$ for all $m$. A
draw from
the final prior $\Pi$ is then constructed as follows:
\begin{enumerate}[(iii)]
\item[(i)] Draw $m$ from the probability distribution $(p_m)$.
\item[(ii)] Independently of $m$, draw a random function from $\Pi'$
and restrict it to $[-m, m]$.
\item[(iii)] Extend the function to the whole real line by setting it
equal to $1$
on $(-\infty, -m-1]$ and to $-1$ on $[m+1, \infty)$, and
interpolating linearly in the intervals $[-m-1,-m]$
and $[m, m+1]$.
\end{enumerate}

By construction, $\Pi$ assigns mass $1$ to a set of drift functions
satisfying the equicontinuity and uniform boundedness
conditions of Theorem \ref{theoremcon}. To prove that this prior yields
consistency for $b_0 \in\mathscr{F}_{s, L}$ it remains
to show that (\ref{eqkl}) holds. Let $\varepsilon> 0$ be fixed. Then
exactly as in the preceding example,
there exists an $m \in\mathbb{N}$ such that
\[
\Pi(b\dvt \|b-b_0\|_{2,\mu_0}^2 < 2\varepsilon^2 )
\ge\Pi(b\dvt \|b-b_0\|_{m,\infty} < \varepsilon).
\]
Since the right-hand side is further bounded from below by
\[
\Pi'(b\dvt \|b - b_0\|_{m, \infty} < \varepsilon) \sum_{n \ge m} p_n,
\]
it now suffices to show that $\Pi'(b\dvt \|b-b_0\|_{m,\infty} <
\varepsilon) > 0$.

To see that this is true,
let $a_k^0$ and $b^0_{j,k}$ be the wavelet coefficients of the
true drift function $b_0$ and let
\[
B = \sum_{k} V_k \varphi_{k} + \sum_{j=0}^J \sum_{k} \eta_j U_{j,k}
\psi_{j,k}
\]
be distributed according to $\Pi'$.
Then
\begin{eqnarray}\label{eqall}
\|B-b_0\|_{m,\infty} &\le&\biggl\| \sum_{k}(V_{k} -a_k^0)
\varphi_k \biggr\|_{m, \infty} \nonumber
\\
 &&{} + \Biggl\| \sum_{j=0}^J \sum_{k } ( \eta_j U_{j,k}
-b_{j,k}^0) \psi_{j,k}\Biggr\|_{m, \infty} +
\biggl\| \sum_{j>J} \sum_{k } b^0_{j,k} \psi_{j,k}\biggr\|_{m, \infty}.
\end{eqnarray}
The first term on the right is bounded by
\[
\|\varphi\|_\infty\sum_{k \in K_m}|V_{k} -a_k^0|,
\]
where $K_m$ is a finite set of natural numbers, since $\varphi$ is
compactly supported.

Since $|a_k^0| \le L$, the $V_k$ have full support in $[-L,L]$ and
$K_m$ is finite, this quantity is bounded by $\varepsilon/3$
with positive probability. By the equivalence of norms mentioned above
and the definition of $\mathscr{F}_{s, L}$,
there exists a constant $c > 0$ such that the
third term on the right of (\ref{eqall}) is bounded by
\[
c \sum_{j>J} 2^{j/2} {\max_{k}} |b^0_{j,k}| \le c\sum_{j>J} 2^{j/2} L
2^{-j(s+1/2)} \le cL 2^{-Js}.
\]
Hence, if we choose $J_0 \in\mathbb{N}$ such that $cL2^{-J_0s} \le
\varepsilon/3$,
then the third term on the right of
(\ref{eqall}) is bounded by $\varepsilon/3$ with probability at
least $\mathbb{P}
(J = J_0) > 0$.
On the event $\{J = J_0\}$,
the second term on the right-hand side of (\ref{eqall}) is bounded
by a constant times
\[
\sum_{j=0}^{J_0} 2^{j/2} \max_{k \in K'_{m}} |\eta_j
U_{j,k}-b^0_{j,k}| \le J_02^{J_0/2} \max_{j \le J_0, k \in K'_m} |\eta
_j U_{j,k}-b^0_{j,k}|.
\]
The set $K'_m$ is finite, since $\psi$ is compactly supported.
Since $|b^0_{j,k}| \le\eta_j L$ and the $U_{j,k}$'s have full support
in $[-L, L]$, the
right-hand side of this display is less than $\varepsilon/3$ as well with
positive probability.
Combining the considerations above and using the fact that $J$, the $V_k$
and the $U_{j,k}$ are all independent, we conclude that
$\Pi'(b\dvt \|b-b_0\|_{m,\infty} < \varepsilon) > 0$.
\end{ex}

\section{\texorpdfstring{Proof of Theorem \protect\ref{theoremcon}}{Proof of Theorem 3.5}}
\label{secproof}

Recall that under $\mathbb{P}_{b}$, the observations $X_0, X_{\Delta
}, \dots
$ form a
discrete-time Markov chain with positive, continuous transition
densities $p_{b}(\Delta, x, y)$ and
a positive, continuous invariant density $\pi_{b}$.
For $b \in\mathscr{B}$, we consider the associated Kullback--Leibler divergence
\[
\mathrm{KL}(b_0, b) = \int\!\!\int p_{b_0}(\Delta, x, y) \log
\frac{p_{b_0}(\Delta, x, y)}{p_{b}(\Delta, x, y)}\pi_{b_0}(x)\, \mathrm{d}x\,\mathrm{d}y.
\]
The following lemma shows that condition (\ref{eqkl}) of Theorem
\ref{theoremcon} implies that we have the Kullback--Leibler
property relative to this distance measure.

\begin{lem}\label{lemkl}
Condition (\ref{eqkl}) of Theorem \ref{theoremcon} implies that for
every $\varepsilon> 0$,
we have $\Pi(b\dvt \operatorname{KL}(b_0, b) < \varepsilon) > 0$.
\end{lem}

\begin{pf}
To prove the lemma we bound
the quantity $\operatorname{KL}(b_0, b)$ from above
by a multiple of $\|b_0-b\|^2_{2, \mu_0}$.
For convenience
we introduce the notation
$K(P,Q) = {\mathbb{E}}_P\log \,\mathrm{d}P/\mathrm{d}Q$ for the Kullback--Leibler divergence
between two probability measures $P$ and $Q$ on the same $\sigma$-field.
The law of a random element $Z$ under the underlying probability
measure $\mathbb{Q}$
is denoted by $\mathscr{L}(Z|\mathbb{Q})$.

Under $\mathbb{P}_b$, for every $b \in\mathscr{B}$, the pair $(X_0,
X_\Delta)$ has
joint density
$(x,y) \mapsto\pi_b(x)p_b(\Delta, x,y)$ relative to Lebesgue measure.
Hence, the Kullback--Leibler divergence between
$\mathscr{L}((X_0,\break
X_\Delta
)|\mathbb{P}_{b_0})$ and
$\mathscr{L}((X_0, X_\Delta)|\mathbb{P}_{b})$ equals
\begin{eqnarray*}
\int\!\!\int\pi_{b_0}(x)p_{b_0}(\Delta, x,y)\log
\frac{\pi_{b_0}(x)p_{b_0}(\Delta, x,y)}{\pi_b(x)p_b(\Delta, x,y)} \,\mathrm{d}x\,\mathrm{d}y
= {\rm KL}(b_0, b) + K(\mu_{b_0}, \mu_b).
\end{eqnarray*}
Now $(X_0, X_\Delta)$ is a measurable functional of the continuous
path $(X_t\dvt t \in[0,\Delta])$. Hence, the Kullback--Leibler divergence
between $\mathscr{L}((X_0, X_\Delta)|\mathbb{P}_{b_0})$ and
$\mathscr{L}((X_0, X_\Delta)|\mathbb{P}_{b})$
is bounded by the Kullback--Leibler divergence between the laws
$\mathscr{L}((X_t\dvt t \in[0,\Delta])|\mathbb{P}_{b_0})$ and
$\mathscr{L}((X_t\dvt t \in[0,\Delta])|\mathbb{P}_{b})$
of the full path $(X_t\dvt t \in[0,\Delta])$
under $\mathbb{P}_{b_0}$ and $\mathbb{P}_b$.
(To see this, observe that the likelihood for $(X_0, X_\Delta)$ is
the conditional expectation of the likelihood for $(X_t\dvt t \in
[0,\Delta])$ and
use the concavity of the logarithm and Jensen's inequality.)
By Girsanov's theorem,
the latter Kullback--Leibler divergence is given by
\[
- {\mathbb{E}}_{b_0}\biggl(\log\frac{\pi_b(X_0)}{\pi_{b_0}(X_0)}
+ \int_0^\Delta(b-b_{0})(X_s) \,\mathrm{d}W_s
- \frac12\int_0^\Delta(b-b_{0})^2(X_s) \,\mathrm{d}s\biggr),
\]
where $W$ is a $\mathbb{P}_{b_0}$-Brownian motion. Using the
stationarity of
the process
$X$ under $\mathbb{P}_{b_0}$, we see that this equals
\[
K(\mu_{b_0}, \mu_b) + \frac{\Delta}{2}\|{b-b_{0}}\|^2_{2, \mu_0}.
\]
Hence, we find that
${2}{\rm KL}(b_0, b) \le{\Delta}\|{b-b_{0}}\|^2_{2, \mu_0}$.
%{\red Misschien de details erbij geven?
\end{pf}

%Let $\scr{F}_{0,\Delta}=\sigma(X_0,X_\Delta)$.
% LHS &= \EE_{b_0} \log\frac{\pi_{b_0}(X_0)p_{b_0}(\Delta, X_0,X_
% \\ &= \EE_{b_0} \log\EE_{b_0} ( \frac{\pi_{b_0}(X_0)}{
%) \\ & \le\EE_{b_0} \EE_{b_0} ( \log\frac{\pi_{b_0}(X_0)}{
%where we applied Jensen's inequality to the log at the inequality.
%Using $d X_t=b_0(X_t) d t + d W_t$ (where $W$ is a $
%& \EE_{b_0}(\log\frac{\pi_b(X_0)}{\pi_{b_0}(X_0)}
%+ \int_0^\Delta(b_0-b)(X_s) dW_s
%+ \frac12\int_0^\Delta(b-b_{0})^2(X_s) ds) \\ & = K(\mu_{b_0},
%using the stationarity of the process
%$X$ under $\PP_{b_0}$.
%}

For any sequence of measurable sets $C_n \subset\mathscr{B}$, we have
that the posterior measure of $C_n$ can be written as
\[
\Pi(C_n|X_0, \ldots, X_{\Delta n}) = \frac{\int_{C_n} L_n(b)
\Pi(\mathrm{d}b)}
{\int_\mathscr{B}L_n(b) \Pi(\mathrm{d}b)},
\]
where
\[
L_n(b) = \frac{\pi_b(X_0)}{\pi_{b_0}(X_0)}\prod_{i=1}^n \frac
{p_b(\Delta, X_{(i-1)\Delta}, X_{i\Delta})}
{p_{b_0}(\Delta, X_{(i-1)\Delta}, X_{i\Delta})}
\]
is the likelihood ratio.
Since we have the Kullback--Leibler property and our Markov chain
satisfies a law of large numbers, the denominator in the expression for
the posterior
can be dealt with in the usual manner. This leads to the following result.

\begin{lem}\label{lemnumerator}
Suppose that for every $\varepsilon> 0$,
we have $\Pi(b\dvt \operatorname{KL}(b_0, b) < \varepsilon) > 0$.
If
for a collection of measurable subsets $C_n \subset\mathscr{B}$ there exists
some $c>0$ such that
\begin{equation}\label{eqpiet}
\mathrm{e}^{nc}\int_{C_n} L_n(b) \Pi(\mathrm{d}b)
\to0,\qquad\mathbb{P}_{b_0}\mbox{-almost surely,}
\end{equation}
then
$\Pi(C_n|X_0, \ldots, X_{\Delta n}) \to0$, $\mathbb
{P}_{b_0}$-almost surely.
\end{lem}

\begin{pf}
By ergodicity, it $\mathbb{P}_{b_0}$-a.s.\ holds that
\begin{eqnarray*}
\frac1n \log L_n(b) = \frac1n\Biggl(\log\frac{\pi_b(X_0)}{\pi_{b_0}(X_0)}+
\sum_{i=1}^n\log\frac{p_b(\Delta, X_{(i-1)\Delta}, X_{i\Delta
})}{p_{b_0}(\Delta, X_{(i-1)\Delta}, X_{i\Delta})}\Biggr)
\to-\operatorname{KL}(b_0, b).
\end{eqnarray*}
In particular, for $\eta> 0$ arbitrary and $b$ such that
$\operatorname{KL}(b_0, b)
< \eta$, it $\mathbb{P}_{b_0}$-a.s. holds that
$\liminf_{n \to\infty} \mathrm{e}^{n\alpha}L_n(b) \ge1$
for all $\alpha> \eta$.
It follows that $\mathbb{P}_{b_0}$-a.s.,
\[
\liminf_{n \to\infty} \mathrm{e}^{n\alpha} {\int_\mathscr{B}L_n(b) \Pi(\mathrm{d}b)}
\ge{\int_{b: \operatorname{KL}(b_0, b) < \eta} \liminf_{n \to\infty}
\mathrm{e}^{n\alpha} L_n(b) \Pi(\mathrm{d}b)}
\ge\Pi\bigl(b\dvt \operatorname{KL}(b_0, b) < \eta\bigr),
\]
and hence
\[
\limsup_{n \to\infty} \Pi(C_n|X_0, \ldots, X_{n\Delta}) \le
\frac{\limsup_{n \to\infty} \mathrm{e}^{{n\alpha}}{\int_{C_n} L_n(b) \Pi(\mathrm{d}b)}}
{\Pi(b\dvt \operatorname{KL}(b_0, b) < \eta)}.
\]
In view of Lemma \ref{lemkl} and the fact that we can take $\alpha>
0$ arbitrarily small, this completes the proof.
\end{pf}

We proceed with the proof of the theorem.
By definition of the topology on $\mathscr{B}$ it suffices to show that
$\Pi(B |X_0, \ldots, X_{n\Delta}) \to0$, $\mathbb{P}_{b_0}$-almost
surely, where\vspace*{2pt}
\[
B = \{b \in\mathscr{B}\dvt \|P^b_\Delta f - P^{b_0}_\Delta f\|_{1, \nu}
> \varepsilon
\},
\]
with $\varepsilon> 0$ and $f$ a continuous function on $\mathbb{R}$
that is
uniformly bounded by $1$.
We fix $\varepsilon$, $f$ and the set $B$ from this point on.

In view of Lemma \ref{lemequi} the assumptions of Theorem \ref{theoremcon},
imply an equicontinuity property for the collections of functions\vspace*{2pt}
\[
\{(P^b_\Delta f)1_K\dvt b\in\mathscr{B}\},
\]
for $K \subset\mathbb{R}$ a compact set.
Arguing as in Tang and Ghosal \cite{Tan07},
this allows us to derive the following useful intermediate result.

\begin{lem}\label{lemhenk}
There exists a compact set $K \subset\mathbb{R}$, a positive integer
$N$ and
bounded intervals $I_{1}, \ldots, I_{N}$ that cover $K$
such that\vspace*{2pt}
\[
B \subset\bigcup_{j=1}^{N} B^+_{j}
\cup\bigcup_{j=1}^{N} B^-_{j},
\]
where\vspace*{2pt}
\begin{eqnarray*}
B^+_{j} & =& \biggl\{b \in B\dvt P^b_\Delta f(x) - P^{b_0}_\Delta f(x) >
\frac\varepsilon{4\nu(K)}\, \forall{x \in I_{j}}\biggr\},
\\[2pt]
B^-_{j} & =& \biggl\{b \in B\dvt P^b_\Delta f(x) - P^{b_0}_\Delta f(x) <-
\frac\varepsilon{4\nu(K)}\, \forall{x \in I_{j}}\biggr\}
\end{eqnarray*}
for $j = 1, \ldots, N$.
\end{lem}

\begin{pf}
Since $\nu$ is a finite Borel measure on the line there exists a
compact subset
$K \subset\mathbb{R}$ such that $\nu(K^c) \le\varepsilon/4$. Let
$\delta>0$ and
cover $K$ with $N< \infty$ intervals with width $\delta/2$,
denote the intervals by $I_{1}, \ldots, I_{N}$. First, we show that $B
\subset\bigcup_{j=1}^{N} B_{j}$,
where\vspace*{2pt}
\[
B_{j} = \biggl\{b \in B\dvt |P^b_\Delta f(x) - P^{b_0}_\Delta f(x)| > \frac
\varepsilon{4\nu(K)} \,\forall{x \in I_{j}}\biggr\}.
\]
Suppose the inclusion is not true. Then there exists a $b\in B$ such
that for each $j\in\{1,\ldots,N\}$ there exists a point $z_j\in I_j$
such that\vspace*{2pt}
\begin{equation}\label{eqzj}
|P^b_\Delta f(z_j) - P^{b_0}_\Delta f(z_j)| \le\frac\varepsilon{4\nu(K)}.
\end{equation}
Now
\begin{eqnarray*}
\|P^b_\Delta f - P^{b_0}_\Delta f\|_{1, \nu}&=&\int_K | P^b_\Delta
f(x)-P^{b_0}_\Delta f(x)| \nu(\mathrm{d}x)+\int_{K^c} | P^b_\Delta
f(x)-P^{b_0}_\Delta f(x)| \nu(\mathrm{d}x)
\\
 &\le&\nu(K) \max_j\max
_{x\in I_j}|P^b_\Delta f(x)-P^{b_0}_\Delta f(x)| + 2\|f\|_\infty\nu
(K^c)
\\
 &\le&
\nu(K)\max_j\max_{x\in I_j} \bigl(|P^b_\Delta f(x)-P^{b}_\Delta
f(z_j)|+|P^b_\Delta f(z_j)-P^{b_0}_\Delta f(z_j)|
\\
&&{}\hspace*{62pt}+|P^{b_0}_\Delta
f(z_j)-P^{b_0}_\Delta f(x)|\bigr)+\varepsilon/2.
\end{eqnarray*}
By local uniform equicontinuity and Lemma \ref{lemequi}, we can find
a $\delta$ such that the first term can be bounded by $\varepsilon
/(8\nu
(K))$. The second term can be bounded by (\ref{eqzj}). By continuity
the third term can be bounded by $\varepsilon/(8\nu(K))$. Therefore, the
preceding display can be bounded by $\varepsilon$, contradicting that
$b\in
B$. Thus $B \subset\bigcup_{j=1}^{N} B_{j}$.

Since the function
$P^b_\Delta f - P^{b_0}_\Delta f$ is continuous and $I_{j}$ is
connected, we have that
$B_{j}$ is included in
\begin{eqnarray*}
&&\biggl\{b \in B\dvt P^b_\Delta f(x) - P^{b_0}_\Delta f(x) > \frac
\varepsilon
{4\nu(K)} \,\forall{x \in I_{j}}\biggr\}
\\
&&\quad{} \cup
\biggl\{b \in B\dvt P^b_\Delta f(x) - P^{b_0}_\Delta f(x) <- \frac
\varepsilon
{4\nu(K)} \,\forall{x \in I_{j}}\biggr\}
=: B_{j}^+ \cup B_{j}^-.
\end{eqnarray*}

%Since $\nu$ is a finite Borel measure on the line there exists a
%compact subset
%$K \subset\RR$ such that $\nu(K^c) \le\eps/4$.
%Then
%for all $b \in B_n$, we have
%and then also
%Let $x_b$ be a point where $|P^b_\Delta f - P^{b_0}_\Delta f|$ attains
%it maximum in $K$.
% By Lemma \ref{lemequi}
%there exists a $\delta> 0$ such that
%for all $b \in B$ and all $x, y \in K$ such that $|x-y| \le\delta$,
%it holds that
%|P^b_\Delta f(x) - P^b_\Delta f(y)| < \frac\eps{8\nu(K)}.
%By the triangle inequality it follows that for all $b \in B$ and for
%all $x \in K$ such that $|x- x_b| < \delta$, we have
%|P^b_\Delta f(x) - P^{b_0}_\Delta f(x)| > \frac\eps{4\nu(K)}.
%Now cover $K$ with $N< \infty$ intervals with width $\delta/2$,
%denote the intervals by $I_{1}, \ldots, I_{N}$.
%Then
%B = \bigcup_{j=1}^{N} \{b \in B: x_b \in I_{j}\} \subset
%where
% B_{j} = \{b \in B: |P^b_\Delta f(x) - P^{b_0}_\Delta f(x)| > \frac
% \]
%Since the function
%$P^b_\Delta f - P^{b_0}_\Delta f$ is continuous and $I_{j}$ is
%connected, we have that
%$B_{j}$ is included in
%& \cup
%=: B_{j}^+ \cup B_{j}^-.
This completes the proof of the lemma.
\end{pf}

As a consequence of this lemma, the proof of the theorem is complete
once we show that for $j = 1, \ldots, N$,
\[
\Pi(B^+_j |X_0, \ldots, X_{n\Delta}) \to0, \qquad
\Pi(B^-_j |X_0, \ldots, X_{n\Delta}) \to0
\]
$\mathbb{P}_{b_0}$-almost surely.
We give the details for the sets $B^+_j$, the argument for the sets
$B^-_j$ is completely analogous. Here, we follow the approach of \cite{Walker}.
We fix $j \in\{1, \ldots, N\}$ and consider the stochastic process
$D$ defined by
\[
D_n = \sqrt{\int_{B_{j}^+} L_n(b) \Pi(\mathrm{d}b)}.
\]
We will show that $\mathbb{P}_{b_0}$-almost surely, $D_n$ converges to $0$
exponentially fast.
According to Lemma~\ref{lemnumerator} this is sufficient.

Note that since $L_n$ is the likelihood, we have ${\mathbb{E}}_{b_0}
D^2_n =
\Pi(B_{j}^+) < \infty$.
Next, we are interested in the conditional expectation ${\mathbb{E}}
_{b_0}(D_{n+1}|\mathscr{F}_n)$,
where $(\mathscr{F}_n)$ is the filtration generated by the Markov chain
$(X_{n\Delta})_{n =0,1, \ldots}$.
Recall that the Hellinger distance $h(p,q)$ between two densities $p,q$
relative to a dominating measure
$\mu$ is
defined by
$h^2(p,q) = \int(\sqrt p-\sqrt q)^2 \,\mathrm{d}\mu$. It satisfies $h^2(p,q) =
2 -2A(p,q)$, where $A(p,q) = \int\sqrt{pq}\, \mathrm{d}\mu$
is the Hellinger affinity between $p$ and $q$.
Then with $p_{n,C}$
the random transition density
\[
p_{n, C}(\Delta, x,y) = \frac{\int_C p_b(\Delta, x, y)L_n(b) \Pi
(\mathrm{d}b)}{\int_C L_n(b) \Pi(\mathrm{d}b)},
\]
we have
\begin{eqnarray*}
{\mathbb{E}}_{b_0}(D_{n+1}|\mathscr{F}_n) & =&
{\mathbb{E}}_{b_0}\Biggl(\sqrt{\int_{B_{j}^+}\frac{p_b(\Delta, X_n,
X_{n+1})}{p_{b_0}(\Delta, X_n, X_{n+1})}
L_n(b) \Pi(\mathrm{d}b)}\Big|\mathscr{F}_n\Biggr)
\\
& =& \int\sqrt{\int_{B_{j}^+}\frac{p_b(\Delta, X_n,
y)}{p_{b_0}(\Delta, X_n, y)}
L_n(b) \Pi(\mathrm{d}b)}p_{b_0}(\Delta, X_n, y) \,\mathrm{d}y
\\
& =& \int\sqrt{\int_{B_{j}^+}{p_b(\Delta, X_n, y)}L_n(b) \Pi
(\mathrm{d}b)p_{b_0}(\Delta, X_n, y)} \,\mathrm{d}y
\\
& =& D_n A_n,
\end{eqnarray*}
where $A_n = A(p_{n,{B_{j}^+}}(\Delta, X_n, \cdot), p_{b_0}(\Delta,
X_n, \cdot))$.
Next, we bound $A_n$. First, note that since $2\|f\|_\infty h(p,q)\ge|
\int f(p-q) \,\mathrm{d} \mu|$, we have $h^2(p,q) \ge\frac14( \int f
(p-q)\, \mathrm{d} \mu)^2$ for functions $f$ that are uniformly bounded by
$1$. Therefore,
\[
A(p,q) =1-\frac12 h^2(p,q) \le1-\frac18\biggl( \int f (p-q)\, \mathrm{d} \mu
\biggr)^2.
\]
Hence, to bound $A_n$ it suffices to lower bound
\[
\int f(y) [ p_{n,B_j^+}(\Delta, X_n,y)-p_{b_0}(\Delta
,X_n,y)]\, \mathrm{d} y
\]
which equals
\[
\int_{B_j^+} \int f(y)[ p_b(\Delta,X_n,y)-p_{b_0}(\Delta
,X_n,y)] \,\mathrm{d}y \frac{L_n(b)}{\int_{B_j^+} L_n(b) \Pi(\mathrm{d}b)} \Pi(\mathrm{d}b).
\]
By the definition of $B_j^+$ in Lemma \ref{lemhenk}, if $X_n\in I_j$
the inner integral is lower bounded by $\varepsilon/(4\nu(K))$. This
implies that
\[
A_n \le1-\frac18\biggl( \frac{\varepsilon}{4\nu(K)} \biggr)^2
\mathbf
{1}_{\{X_n \in I_j\}}.
\]
%
%It follows immediately from the definitions of ${B_{j}^+}$ and
%$p_{n,C}$ that for $x \in I_j$,
%Since for densities $p,q$ it holds that
%$|\int f(p-q) d\mu| \le2\|f\|_\infty h(p,q)$,
%it follows that for $x \in I_j$,
%we have $h(p_{n,{B_{j}^+}}(x, \cdot), p_{b_0}(\Delta, x, \cdot)) \ge
%But then
%$A_n \le1-(\eps^2/8)1_{X_n \in I_j}$
Hence,
\[
{\mathbb{E}}_{b_0}(D_{n+1}|\mathscr{F}_n) \le D_n
(1-k\varepsilon^2 1_{X_n \in
I_j}),
\]
where $k=1/(128 \nu(K)^2)$.
We conclude that the process
\[
M_n = D_n(1-k\varepsilon^2)^{-\sum_{i=1}^{n-1} 1_{X_i \in I_j}}
\]
is an $(\mathscr{F}_n)$-supermartingale under the measure $\mathbb
{P}_{b_0}$ (note that
$M_n$ is bounded by the integrable process $D_n(1-k\varepsilon^2)^{-(n-1)}$,
hence $M_n$ is integrable).
By Doob's martingale convergence theorem, we have $M_n \to M_\infty$
almost surely,
for some finite-valued random variable $M_\infty$.
By ergodicity, we have $n^{-1}\sum_{i=1}^{n-1} 1_{X_i \in I_j} \to\mu
_{b_0}(I_j) > 0$
almost surely.
An application of Lemma \ref{lemnumerator} completes the proof.

\section{Concluding remarks}\label{sec6}

In this paper, we obtain conditions for posterior consistency of nonparametric
Bayesian drift estimation for low-frequency observations from a scalar
ergodic diffusion.
The main theorem and the subsequent examples provide several priors for
which consistency
is guaranteed.
As discussed in the \hyperref[sec1]{Introduction}, data augmentation techniques that
have been proven
to be effective in parametric settings, are in principle usable for
numerical implementation of nonparametric models as well.
Preliminary investigations indicate that practically feasible
procedures can indeed be constructed, but more
work on computational issues is necessary at the moment.

The results and proofs in this paper show that in this low-frequency
observations setting, obtaining consistency
relative to a rather weak topology is already quite involved.
Very challenging but equally interesting would be the development of a
testing approach to posterior consistency in this
setting. It would allow to obtain consistency in stronger topologies,
rates of contraction and
relaxation of boundedness and equicontinuity conditions. For general
diffusions this seems rather difficult, but
some progress might be possible for diffusions on compact state spaces.

\begin{appendix}

\section*{Appendix: An equicontinuity property of the transition operators}\label{app}

The concept of local uniform equicontinuity is given in Definition \ref{dflocunifequicont}.
\begin{lem}\label{lemequi}
If $\sup_{b \in\mathscr{B}} \|b\|_\infty< \infty$ and
$\mathscr{B}$ is locally uniformly equicontinuous,
then for every $f \in BC(\mathbb{R})$ and $t > 0$, the collection
$\{P^b_t f\dvt b \in\mathscr{B}\}$ is locally uniformly equicontinuous
as well.
\end{lem}

\begin{pf}
Let $K \subset\mathbb{R}$ be a compact set.
For $\mathbb{P}^0_x$ the law of the Brownian motion
starting in $x$ we have, by Girsanov's theorem,
\begin{eqnarray*}
P^b_tf(x) = {\mathbb{E}}_x^0 f(X_t)\frac{\mathrm{d}\mathbb{P}^b_x}{\mathrm{d}\mathbb
{P}^0_x} = {\mathbb{E}}_x^0 f(X_t)
\exp\biggl(\int_0^{t}b(X_s)\, \mathrm{d}X_s-\frac12 \int_0^{t}b^2(X_s)\, \mathrm{d}s
\biggr).
\end{eqnarray*}
Under $\mathbb{P}^0_x$ the process $X$ has the same law as $x+W$, for
$W$ a standard
Brownian motion starting in $0$. Hence, we get
\[
P^b_tf(x) = {\mathbb{E}}f(x +W_t) L_x,
\]
where
\[
L_u = \mathrm{e}^{l_u}, \qquad l_u = \int_0^{t}b(u+W_s) \,\mathrm{d}W_s-\frac12 \int
_0^{t}b^2(u+W_s)\, \mathrm{d}s.
\]
It follows that
\begin{eqnarray*}
|P^b_t f(x) - P^b_tf(y)| & \le&{\mathbb{E}}|f(x +W_t) L_x -f(y +W_t)
L_y|
\\
& \le&{\mathbb{E}}|f(x +W_t)| |L_x-L_y| +{\mathbb{E}}|L_y||f(x +W_t)
- f(y +W_t)|
\\
&\phantom{.} =:& I + \mathit{II}.
\end{eqnarray*}

We first bound the term $I$. By the fact that $|\mathrm{e}^b-\mathrm{e}^a| \le
|a-b|(\mathrm{e}^a+\mathrm{e}^b)$ and Cauchy--Schwarz,
\begin{eqnarray*}
|I|^2 & \le&\|f\|^2_\infty({\mathbb{E}}|L_x-L_y|)^2
 \\
& \le&
\|f\|^2_\infty({\mathbb{E}}|l_x-l_y||L_x+L_y|)^2
 \\
& \le&\|f\|^2_\infty{\mathbb{E}}|l_x-l_y|^2 {\mathbb{E}}(L_x + L_y)^2.
\end{eqnarray*}
We have
\begin{equation}\label{eql}
l_x - l_y = \int_0^{t}\bigl(b(x+W_s)-b(y+W_s)\bigr)\, \mathrm{d}W_s-\frac12 \int
_0^{t}\bigl(b^2(x+W_s)-b^2(y+W_s)\bigr)\, \mathrm{d}s.
\end{equation}
For the first term on the right the It\^o isometry gives, for $x, y \in K$,
\begin{eqnarray*}
&&{\mathbb{E}} \biggl(\int_0^{t}\bigl(b(x+W_s)-b(y+W_s)\bigr)\, \mathrm{d}W_s\biggr)^2
\\
&&\quad =
{\mathbb{E}}\int_0^{t}\bigl(b(x+W_s)-b(y+W_s)\bigr)^2 \,\mathrm{d}s
\\
&&\quad  \le{\mathbb{E}}\int_0^{t}\bigl(b(x+W_s)-b(y+W_s)\bigr)^21_{\sup_{s \le t}
|W_s| \le
M} \,\mathrm{d}s
+ 4t\|b\|^2_\infty\mathbb{P}\Bigl({\sup_{s \le t} |W_s| > M}\Bigr)
\\
&&\quad  \le{t}\mathop{\sup_{u,v \in K' }}_{|u-v| \le|x-y|}|b(u)-b(v)|^2
+ 4t\|b\|^2_\infty\mathbb{P}\Bigl({\sup_{s \le t} |W_s| > M}\Bigr)
\end{eqnarray*}
for every $M > 0$, where $K' =\{x +y\dvt x \in K, y \in[-M,M]\}$.
The assumptions on $\mathscr{B}$ imply that by choosing $M$ large
enough and
$|x-y|$ small enough,
the right-hand side can be made arbitrarily small, uniformly in
$\mathscr{B}$.
The second term on the right of (\ref{eql}) can be handled in the
same manner,
using also the fact that $|b^2(u)-b^2(v)| \le2\|b\|_\infty|b(u)-b(v)|$.
To complete the bound for term $I$, we
note that ${\mathbb{E}}(L_x + L_y)^2 \le2{\mathbb{E}}L_x^2 +
2{\mathbb{E}}L^2_y$ and we write
\begin{eqnarray*}
L_u^2 & =& \exp\biggl({\int_0^{t}2b(u+W_s) \,\mathrm{d}W_s- \int
_0^{t}b^2(u+W_s)\, \mathrm{d}s}\biggr)
\\
& =& \exp\biggl({\int_0^{t}b^2(u+W_s)\, \mathrm{d}s}\biggr) \exp\biggl({\int
_0^{t}2b(u+W_s)\, \mathrm{d}W_s- \frac12\int_0^{t}(2b)^2(u+W_s)\, \mathrm{d}s}
\biggr).
\end{eqnarray*}
The first factor on the right is bounded by $\exp(t\|b\|^2_\infty)$
and the second one
is the time $t$ value of a martingale that starts in $1$. Hence,
\[
{\mathbb{E}}L_u^2 \le \mathrm{e}^{t\|b\|^2_\infty}.
\]
Finally, observe that
by Cauchy--Schwarz and a bound derived above,
\[
|\mathit{II}|^2 \le \mathrm{e}^{t\|b\|^2_\infty}{\mathbb{E}}|f(x +W_t) - f(y +W_t)|^2.
\]
This completes the proof.
\end{pf}

\end{appendix}

\section*{Acknowledgements} We thank one of the referees for helpful
comments on Example 4.2. The second author's research is partially funded by the
Netherlands Organization for Scientific Research (NWO).

% imsref loaded by elazauskaite, 2011-11-21 09:09:59
%

\printhistory

\end{document}